\documentclass[12pt]{article}
\usepackage{amsbsy,amsfonts,amssymb}
\usepackage{graphicx}

\newtheorem{theo}{Theorem}[section]
\newtheorem{lemma}[theo]{Lemma}

\newtheorem{example}[theo]{Example}
\newtheorem{defi}[theo]{Definition}
\newtheorem{coll}[theo]{Corollary}
\newtheorem{remark}[theo]{Remark}

\newcommand{\be}{\begin{equation}}
\newcommand{\ee}{\end{equation}}
\newcommand{\ba}{\begin{array}}
\newcommand{\ea}{\end{array}}

\newcommand{\dsum}{\displaystyle \sum}
\newcommand{\dprod}{\displaystyle \prod}

\begin{document}

\title{Coloring of graphs associated to zero-divisors}

\author{
Hsin-Ju Wang 
  }

\date{}

\maketitle

\noindent{\bf Abstract.} Let $G$ be a graph, $\chi(G)$ be the
minimal number of colors which can be assigned to the vertices of
$G$ in such a way that every two adjacent vertices have different
colors and $\omega(G)$ to be the least upper bound of the size of
the complete subgraphs contained in $G$. It is well-known that
$\chi(G)\geq \omega(G)$. Beck in \cite{b} conjectured that
$\chi(\Gamma_0(R))=\omega(\Gamma_0(R))$ if
$\omega(\Gamma_0(R))<\infty$, where $\Gamma_0(R)$ is a graph
associated to a commutative ring $R$. In this note, we provide
some sufficient conditions for a ring $R$ to enjoy
$\chi(\Gamma_0(R))=\omega(\Gamma_0(R))$. As a consequence, we
verify Beck's conjecture for the homomorphic image of
$\mathbb{Z}^n$.

\section{Introduction}
In this paper, $R$ will denote a commutative ring with $1$. If $S$
is a subset of $R$, we denote $S-\{0\}$ by $S^*$. Also, we use
$\mathbb{N}$ for the set of all nonnegative integers. \par In
\cite{b}, Beck introduced the idea of a {\it zero-divisor graph}
of a commutative ring $R$ with 1. He defined $\Gamma_0(R)$ to be
the graph whose vertices are elements of $R$ and in which two
vertices $x$ and $y$ are adjacent if and only if $xy=0$. Beck was
mostly concerned with coloring $\Gamma_0(R)$. Recall that the {\it
chromatic number} of the graph, denoted $\chi(\Gamma_0(R))$, is
defined to be the minimal number of colors which can be assigned
to the elements of $\Gamma_0(R)$ in such a way that every two
adjacent elements have different colors. A subset $C$ of vertex
set of $\Gamma_0(R)$ is called a {\it clique} if any two distinct
elements of $C$ are adjacent; the number $\omega(\Gamma_0(R))$ is
the least upper bound of the size of the cliques, and clearly
$\chi(\Gamma_0(R))\geq \omega(\Gamma_0(R))$. Beck conjectured that
$\chi(\Gamma_0(R))=\omega(\Gamma_0(R))$ if
$\omega(\Gamma_0(R))<\infty$ and proved this conjecture for a
rather wide class of rings, including reduced rings and principal
ideal rings. However, D.D. Anderson and M. Naseer gave a
counterexample to this conjecture in \cite{an}. One of the aim of
this paper is to find some sufficient conditions for rings $R$ to
enjoy $\chi(\Gamma_0(R))=\omega(\Gamma_0(R))$. As a consequence,
we are able to verify the conjecture for a class of rings, namely,
the homomorphic image of $\mathbb{Z}^n$.
\par Let $R$ be a ring and let $Z(R)$ be the
set of zero-divisors of $R$. Then there is associated a (simple)
graph $\Gamma(R)$ to $R$ with vertices $Z(R)^*$ and where two
distinct vertices $x,y\in Z(R)^*$ are adjacent if and only if
$xy=0$. This graph was defined slightly differently than the graph
introduced by I. Beck and was introduced by Anderson and
Livingston in \cite{al}. Recently, this concept was studied
extensively in \cite{amy}, \cite{am}, \cite{afll}, \cite{al},
\cite{s} and \cite{wa}. One of the aim of this paper is to  prove
that $\chi(\overline{\Gamma(R)})=\omega(\overline{\Gamma(R)})$ if
$R\cong \mathbb{Z}_{p_1^{r_1}}\times \cdots \times
\mathbb{Z}_{p_k^{r_k}}$, where $p_i$ is a prime number for every
$i$ and $\bar{G}$ is the complement graph of $G$.

\section{Preliminary}
We review some background from graph theory and fix some notations
from \cite{an} and \cite{b} in this section. \par A simple graph
$G$ is an ordered pair of disjoint sets $(V, E)$ such that
$V=V(G)$ is the vertex set of $G$ and $E=E(G)$ is the edge set of
$G$. Often we use $G$ for $V(G)$. The {\it order} of a graph $G$,
written by $|G|$, is the cardinality of $V(G)$. A {\it subgraph}
of $G$ is a graph having all of its vertices and edges in $G$. Let
$V'\subseteq V(G)$; then $G-V'$ is the subgraph of $G$ obtained by
deleting the vertices in $V'$ and all edges incident with them.
For $v\in V$, the {\it degree} of $v$, denoted by $deg(v)$, is the
number of edges of $G$ incident to $v$. A vertex subset $S$ of a
graph $G$ is called an {\it independent set} if any two vertices
of $S$ are not adjacent. \par A graph $G$ is said to be a {\it
totally disconnected} graph if the edge set $E(G)=\emptyset$. A
graph in which each pair of distinct vertices is joined by an edge
is called a {\it complete graph}. We use $K_n$ for the complete
graph with $n$ vertices.
\par The most important result in \cite{b}
is the following.

\begin{theo}
\label{btheo} \cite[Theorem~3.9]{b} The following are equivalent
for a ring $R$:
\begin{description}
\item{(i)} $\chi(\Gamma_0(R))<\infty$.  \item{(ii)}
$\omega(\Gamma_0(R))<\infty$. \item{(iii)} The nil-radical is
finite and equals a finite intersection of prime ideals.
\end{description}
\end{theo}

Beck named the rings satisfying one of the three conditions in
Theorem\ref{btheo} as follows.
\begin{defi} \label{defi1}
A ring $R$ is called a coloring provided that
$\chi(\Gamma_0(R))<\infty$.
\end{defi}
In case $\chi(\Gamma_0(R))=\omega(\Gamma_0(R))<\infty$ for rings
$R$, Anderson and Naseer made the following definition.
\begin{defi} \label{defi2}
A ring $R$ is called a chromatic provided that
$\chi(\Gamma_0(R))=\omega(\Gamma_0(R))<\infty$.
\end{defi}
\par In the sequel, we use the symbol $A\sqcup B$ to denote
the {\it disjoint union} of two sets $A$ and $B$.

\section{Coloring of $\Gamma_0(R)$}

Let $R$ be a commutative ring with 1 and let $Z(R)$ be the set of
zero-divisors of $R$. Let $\Gamma_0(R)$ and $\Gamma(R)$ be the
graphs associated to $R$ defined in the introduction. Observe
that, in $\Gamma_0(R)$, $0$ is adjacent to every other vertices of
$\Gamma_0(R)$ and the degree of every unit is $1$. Thus
$\Gamma(R)$ is a subgraph of $\Gamma_0(R)$ and $\Gamma_0(R)-\{0\}$
is a disjoint union of $\Gamma(R)$ and a totally disconnected
graph. \par  The goal of this section is to study the interplay
between $\chi(\Gamma_0(R))$ and $\omega(\Gamma_0(R))$. For this,
we first look at some special cases.

\begin{lemma}
\label{lem1} Let $R\cong \mathbb{Z}_{p^r}$, where $p$ is a prime
number and $r$ be a positive integer. Then $R$ is a chromatic ring
with $\chi(\Gamma_0(R))=\omega(\Gamma_0(R))=p^{r/2}$if $r$ is even
and $\chi(\Gamma_0(R))=\omega(\Gamma_0(R))=p^{(r-1)/2}+1$ if $r$
is odd.
\end{lemma}
\begin{proof}
Suppose that $r=2t$ is even. Then $S=\{p^tn~|~n\in \mathbb{Z}\}$
is a clique subset of $\Gamma_0(R)$. Since $|S|=p^{r/2}$, we
require $p^{r/2}$ colors to paint $S$. Observe that $p^t$ is not
adjacent to any vertex of $\Gamma_0(R)-S$ and $\Gamma_0(R)-S$ is
an independent set. We can use the color of $p^t$ to color every
vertex of $\Gamma_0(R)-S$. Thus $\chi(\Gamma_0(R))\leq p^{r/2}\leq
\omega(\Gamma_0(R))\leq \chi(\Gamma_0(R))$.\par Suppose that
$r=2t+1$ is odd. Then $S=\{p^{t+1}n~|~n\in \mathbb{Z}\}\cup
\{p^t\}$ is a clique subset of $\Gamma_0(R)$. Since
$|S|=p^{(r-1)/2}+1$, we require $p^{(r-1)/2}+1$ colors to paint
$S$. Observe that $p^t$ is not adjacent to any vertex of
$\Gamma_0(R)-S$ and $\Gamma_0(R)-S$ is an independent set. We can
use the color of $p^t$ to color every vertex of $\Gamma_0(R)-S$.
Thus $\chi(\Gamma_0(R))\leq p^{(r-1)/2}+1\leq
\omega(\Gamma_0(R))\leq \chi(\Gamma_0(R))$.
\end{proof}

From the proof of Lemma~\ref{lem1}, we have the following
consequence.

\begin{coll} \label{coll1}
Let $R\cong \mathbb{Z}_{p^r}$, where $p$ is a prime number and $r$
be a positive integer. Let $S$ be a largest clique of
$\Gamma_0(R)$; then $|\{a\in S~|~a^2\neq 0\}|=0$ if $r$ is even
and $|\{a\in S~|~a^2\neq 0\}|=1$ if $r$ is odd.
\end{coll}

\begin{example}
\label{exam1} Let $R\cong \mathbb{Z}_8\times \mathbb{Z}_{16}$;
then $\chi(\Gamma_0(R))=\omega(\Gamma_0(R))=9$.
\end{example}
\begin{proof}
Let $C_1=\{(0, 0), (0, 4), (0, 8), (0, 12), (4, 0), (4, 4), (4,
8), (4, 12), (2, 0)\}$, $C_2=(\{1, 3, 5, 7\}\times \{0, 4, 8,
12\})\cup ((\{2, 6\}\times \{0, 4, 8, 12\})-\{(2, 0)\})$ and
$C_3=\mathbb{Z}_8\times (\mathbb{Z}_{16}-\{0, 4, 8, 12\})$; then
$C_1$ is a clique and $\Gamma_0(R)=C_1\sqcup C_2\sqcup C_3$.
Moreover, $C_2$ and $C_3$ are independent sets. Observe that we
require $9$ colors to paint $C_1$. Then we can use the color of
$(2, 0)$ (resp. $(0, 4)$)  to color $C_2$ (resp. $C_3$). In such a
way, every two adjacent vertices of $\Gamma_0(R)$ have different
colors. Thus, we conclude that $9\leq \omega(\Gamma_0(R))\leq
\chi(\Gamma_0(R))\leq 9$.
\end{proof}

The following result stated in \cite{b} without any proof. Here,
we provide one.
\begin{lemma}
\label{lem2} Let $R_1$ and $R_2$ be rings and $R\cong R_1\times R_2$. Then the following hold. \\
(i) $\chi(\Gamma_0(R))\geq
\chi(\Gamma_0(R_1))+\chi(\Gamma_0(R_2))-1$.\\ (ii) If $R_2$ is
reduced, then
$\chi(\Gamma_0(R))=\chi(\Gamma_0(R_1))+\chi(\Gamma_0(R_2))-1$.
\end{lemma}
\begin{proof}
(i) If  $\chi(\Gamma_0(R_1))$ or $\chi(\Gamma_0(R_2))$ are
infinite then so is $\chi(\Gamma_0(R))$. Therefore the inequality
holds trivially. So we may assume that
$t_i=\chi(\Gamma_0(R_i))<\infty$ for $i=1, 2$. Thus, we require
$t_1$ colors to paint $\Gamma_0(R_1)\times \{0\}$ and $t_2$ colors
to paint $\{0\}\times \Gamma_0(R_2)$. Since $(0, 0)$ is adjacent
to any vertex of $\Gamma_0(R)$ and every vertex of
$\Gamma_0(R_1)\times \{0\}$ is adjacent to every vertex of
$\{0\}\times \Gamma_0(R_2)$, we require at least $t_1+t_2-1$
colors to paint $\Gamma_0(R)$.\\ (ii) Suppose that $R_2$ is
reduced.  As in (i), we may assume that
$t_i=\chi(\Gamma_0(R_i))<\infty$ for $i=1, 2$. Observe that
$t_i-1=\chi(\Gamma_0(R_i)-\{0\})$ for $i=1, 2$, there are
independent sets $S_1, \dots, S_{t_1-1}$ of the graph
$\Gamma_0(R_1)$  and independent sets $T_1, \dots, T_{t_2-1}$ of
the graph $\Gamma_0(R_2)$ such that $\Gamma_0(R_1)-\{0\}=S_1\sqcup
\cdots \sqcup S_{t_1-1}$ and $\Gamma_0(R_2)-\{0\}=T_1\sqcup \cdots
\sqcup T_{t_2-1}$. Thus, \be \label{eq1} \ba{rl} \Gamma_0(R)= &
(0,
0)\sqcup (S_1\times \{0\})\sqcup \cdots (S_{t_1-1}\times \{0\}) \\
& \sqcup (\Gamma_0(R_1)\times T_1)\sqcup \cdots \sqcup
(\Gamma_0(R_1)\times T_{t_2-1}). \ea \ee Observe that
$\Gamma_0(R_1)\times T_i$ is an independent set in $\Gamma_0(R)$
for every $i$: Let $(a_j, b_j)\in \Gamma_0(R_1)\times T_i$ for
$j=1, 2$. If $(a_1, b_1)(a_2, b_2)=(0, 0)$, then $b_1b_2=0$, so
that $b_1^2=0$ if $b_1=b_2$ or $b_1$ is adjacent to $b_2$ in $T_i$
if $b_1\neq b_2$. Both statements lead to contradiction.
\par Notice that by (\ref{eq1}), $\Gamma_0(R)$ is a disjoint union
of $t_1+t_2-1$ independent sets. Now, we use $t_1+t_2-1$ colors to
color those independent sets. In this way, one can easily to check
that every two adjacent vertices have different colors. Thus,
$\chi(\Gamma_0(R))\leq \chi(\Gamma_0(R_1))+\chi(\Gamma_0(R_2))-1$.
\end{proof}

Similar to Lemma~\ref{lem2}, we have the following.
\begin{lemma}
\label{lem3} Let $R_1$ and $R_2$ be rings and $R\cong R_1\times R_2$. Then the following hold. \\
(i) $\omega(\Gamma_0(R))\geq
\omega(\Gamma_0(R_1))+\omega(\Gamma_0(R_2))-1$.\\ (ii) If $R_2$ is
reduced, then
$\omega(\Gamma_0(R))=\omega(\Gamma_0(R_1))+\omega(\Gamma_0(R_2))-1$.
\end{lemma}
\begin{proof}
(i) If  $\omega(\Gamma_0(R_1))$ or $\omega(\Gamma_0(R_2))$ are
infinite then so is $\omega(\Gamma_0(R))$. Therefore the
inequality holds trivially. So we may assume that
$t_i=\omega(\Gamma_0(R_i))<\infty$ for $i=1, 2$. Therefore, there
is a maximal clique $S$ (resp. $T$) in $\Gamma_0(R_1)$ (resp.
$\Gamma_0(R_2)$) such that $|S|=t_1$ (resp. $|T|=t_2$). Notice
that $\{(0, 0)\}\cup ((S-\{0\})\times \{0\})\cup (\{0\}\times
(T-\{0\}))$ is a clique in $\Gamma_0(R)$. Thus,
$\omega(\Gamma_0(R))\geq
\omega(\Gamma_0(R_1))+\omega(\Gamma_0(R_2))-1$.\\ (ii) Suppose
that $R_2$ is reduced. As in (i), we may assume that
$t_i=\omega(\Gamma_0(R_i))<\infty$ for $i=1, 2$. Let $W$ be any
maximal clique of $\Gamma_0(R)$. Let $s=|W\cap
(\Gamma_0(R_1)\times (\Gamma_0(R_2)-\{0\}))|$ and $\{(a_1, b_1),
\dots, (a_s, b_s)\}=W\cap (\Gamma_0(R_1)\times
(\Gamma_0(R_2)-\{0\}))$. Notice that $b_ib_j=0$ for every $i\neq
j$. If $b_i=b_j$ for some $i\neq j$, then $b_i^2=0$, so that
$b_i=0$ as $R_2$ is reduced, a contradiction. Therefore $\{0, b_1,
\dots, b_s\}$ is a clique of $\Gamma_0(R_2)$. It follows that
$s+1\leq t_2$. On the other hand, $W\cap (\Gamma_0(R_1)\times
\{0\})$ is a clique of $\Gamma_0(R_1)\times \{0\}$. Thus, $|W\cap
(\Gamma_0(R_1)\times \{0\})|\leq t_1$. Hence we conclude that
$|W|\leq t_1+t_2-1$ and then $\omega(\Gamma_0(R))\leq
\omega(\Gamma_0(R_1))+\omega(\Gamma_0(R_2))-1$.
\end{proof}

An immediately consequence is the following.
\begin{coll}
\label{coll1} Let $R\cong R_1\times \cdots \times R_k$, where
$R_i$ is an integral domain for every $i$; then
$\chi(\Gamma_0(R))=\omega(\Gamma_0(R))=k+1$.
\end{coll}

We now state and prove our main result in this section.

\begin{theo}
\label{maintheo1} Let $R\cong R_1\times \cdots \times R_k$, where
$R_i$ is a coloring  ring for every $i$. Suppose that there is a
finite set $S_i$ of $\Gamma_0(R_i)$ for every $i$ such that the
following hold:
\begin{description}
\item{(i)} $S_i$ is a maximal clique of $\Gamma_0(R_i)$.
\item{(ii)} If $a, b\in \Gamma_0(R_i)-S_i$, then $a^2\neq 0$ and
$ab\neq 0$.
\end{description}
Let $N_i=\{a\in S_i~|~a^2=0\}\subseteq S_i$ and $n_i=|N_i|$ for
every $i$; then
$\chi(\Gamma_0(R))=\omega(\Gamma_0(R))=\dprod_{i=1}^k
n_i+\dsum_{i=1}^k (s_i-n_i)$, where $s_i=|S_i|$.
\end{theo}
\begin{proof}
Rearrange $R_i$ if necessary, we may assume that $S_i=N_i$ if and
only if $i>t$ for some nonnegative integer $t\leq k$. For $i=1,
\dots, t$, let $$A_i=\{(a_1, \dots, a_k)\in
\Gamma_0(R)~|~a_j=0~if~j\neq i~and~a_i\in S_i-N_i\}$$ and
$A=\cup_{i=1}^t A_i$. Let $$N=N_1\times \cdots \times N_k,$$
$$S=S_1\times \cdots \times S_k,$$
$$T_1=(Z(R_1)\times \cdots \times Z(R_k))-S$$ and $$T_2=(\Gamma_0(R_1)\times
\cdots \times \Gamma_0(R_k))-T_1;$$ then $\Gamma_0(R)=S\sqcup
T_1\sqcup T_2$. Further, $N\cup A$ is a clique. Let
$n=\prod_{i=1}^k n_i+\sum_{i=1}^k (s_i-n_i)$; then $|N\cup A|=n$
and we require $n$ colors to paint $N\cup A$ as $N\cup A$ is a
clique. We will show that $\chi(\Gamma_0(R))\leq n$. For this, we
fix some notations: \par  For every $i=1, \dots, k$, choose
$b_i\in S_i-\{0\}$. This can be done as $|S_i|\geq 2$.
\par For every $i=1, \dots, k$ and $c\in Z(R_i)-S_i$, there are
vertices in $S_i$ that are not adjacent to $c$ by assumption (i).
Choose an element enjoys this property and denote this element as
$\tilde{c}$. So $c\tilde{c}\neq 0$ and $\tilde{c}\in S_i$. \par If
$(u_1, \dots, u_k)\in N\cup A$ and $T$ is an independent set of
$\Gamma_0(R)$, then we use $T\thicksim \{(u_1, \dots, u_k)\}$ to
denote the fact that every vertex of $T$ is colored by the color
of $(u_1, \dots, u_k)$. Furthermore, if $U$ and $V$ are two
vertices of $\Gamma_0(R)$, then we use $U\thicksim V$ (resp.
$U\nsim V$) to denote the fact that the colors of $U$ and $V$ are
the same (resp. different).
\par To finish the proof, we need to color
every vertex of $S$, $T_1$ and $T_2$.\\ \underline{Coloring of
$S$}: \par First, we use $n$ colors to color $N\cup A$. Notice
that
$$S-N=\sqcup_{i=1}^t \sqcup_{a\in S_i-N_i} N_1\times \cdots \times
N_{i-1}\times \{a\}\times S_{i+1}\times \cdots \times S_k.$$
Notice also that $N_1\times \cdots \times N_{i-1}\times
\{a\}\times S_{i+1}\times \cdots \times S_k$ is an independent set
for every $a\in S_i-N_i$ and for every $i$. Now, we can use the
color of $(0, \cdots, 0, a, 0 \cdots, 0)$ to color the above
independent set, that is, \be \label{eq2} N_1\times \cdots \times
N_{i-1}\times \{a\}\times S_{i+1}\times \cdots \times S_k\thicksim
\{(0, \cdots, 0, a, 0 \cdots, 0)\}. \ee ~~\\ \underline{Coloring
of $T_1$}: \par Observe that $$T_1=\sqcup_{i=1}^k \sqcup_{c\in
Z(R_i)-S_i} S_1\times \cdots \times S_{i-1}\times \{c\}\times
Z(R_{i+1})\times \cdots \times  Z(R_k).$$ Moreover, $S_1\times
\cdots \times S_{i-1}\times \{c\}\times Z(R_{i+1})\times \cdots
\times   Z(R_k)$ is an independent set for every $c\in Z(R_i)-S_i$
and for every $i$. Now, we can use the color of $(0, \cdots, 0,
\tilde{c}, 0 \cdots, 0)$ to color the above independent set, that
is, \be \label{eq3} S_1\times \cdots \times S_{i-1}\times
\{c\}\times
Z(R_{i+1})\times \cdots \times Z(R_k)\thicksim \{(0, \cdots, 0, \tilde{c}, 0 \cdots, 0)\}. \ee ~~\\
\underline{Coloring of $T_2$}: \par Observe that
$$T_2=\sqcup_{i=1}^k \sqcup_{b\in \Gamma_0(R_i)-Z(R_i)} Z(R_1)\times \cdots
\times Z(R_{i-1})\times \{b\}\times \Gamma_0(R_{i+1})\times \cdots
\times \Gamma_0(R_k).$$ Moreover, $Z(R_1)\times \cdots \times
Z(R_{i-1})\times \{b\}\times \Gamma_0(R_{i+1})\times \cdots \times
\Gamma_0(R_k)$ is an independent set for every $b\in
\Gamma_0(R_i)-Z(R_i)$ and for every $i$. Now, we can use the color
of $(0, \cdots, 0, b_i, 0 \cdots, 0)$ to color the above
independent set, that is, \be \label{eq4} Z(R_1)\times \cdots
\times Z(R_{i-1})\times \{b\}\times \Gamma_0(R_{i+1})\times \cdots
\times \Gamma_0(R_k))\thicksim \{(0, \cdots, 0, b_i, 0 \cdots,
0)\}. \ee \par To complete the proof, we need to verify that every
two adjacent vertices in $\Gamma_0(R)$ have different colors. For
this, let $C=(c_1, \dots, c_k)$ and $D=(d_1, \dots, d_k)$ be two
adjacent vertices in $\Gamma_0(R)$; then $c_id_i=0$ for every $i$.
We proceed by discussing the following cases. \\ Case~1. $C, D\in
S$. Assume first that $C, D\in N$. Then $C\nsim D$ of course.
Assume next that $C\in N$ and $D\notin N$. Let $i$ be the least
integer for which $d_i\notin N_i$; then $D\thicksim D'=(0, \cdots,
0, d_i, 0 \cdots, 0)$ by (\ref{eq2}). Since $c_id_i=0$, it follows
that $C$ and $D'$ are different vertices in $N\cup A$. Therefore,
$C\nsim D'$ and then $C\nsim D$. Finally, assume that $C\notin N$
and $D\notin N$. Let $i$ be the least integer for which $c_i\notin
N_i$ and let $j$ be the least integer for which $d_j\notin N_j$;
then $C\thicksim C'=(0, \cdots, 0, c_i, 0 \cdots, 0)$ and
$D\thicksim D'=(0, \cdots, 0, d_j, 0 \cdots, 0)$ by (\ref{eq2}).
If $i=j$, then $c_id_i=0$ implies that $C'\nsim D'$. Thus, $C\nsim
D$. If $i\neq j$, then of course $C'$ and $D'$ are different
vertices in $N\cup A$, so that $C'\nsim D'$, it follows that
$C\nsim D$.\\ Case~2. $C\in S$ and $D\in T_1$. Let $i$ be the
least integer for which $d_i\notin S_i$; then $D\thicksim D'=(0,
\cdots, 0, \tilde{d_i}, 0 \cdots, 0)$ by (\ref{eq3}). Furthermore,
$c_id_i=0$ and $d_i\tilde{d_i}\neq 0$ implies that $c_i$ and
$\tilde{d_i}$ are different vertices in $S_i$. Since $S_i$ is a
clique, $c_i\tilde{d_i}=0$. If $C\in N$, then $C\nsim D'$. If
$C\notin N$, then $C\thicksim C'= (0, \cdots, 0, c_j, 0 \cdots,
0)$ by
(\ref{eq2}), so that $C'$ and $D'$ are different vertices in $N\cup A$, it follows that $C\nsim D$.\\
Case~3. $C\in S$ and $D\in T_2$. Let $i$ be the least integer for
which $d_i\notin Z(R_i)$; then $D\thicksim D'=(0, \cdots, 0, b_i,
0 \cdots, 0)$ by (\ref{eq4}). Furthermore, $c_id_i=0$ and $d_i$ is
a non-zero-divisor implies that $c_i=0$ and then $c_i$ and $b_i$
are different vertices in $S_i$. If $C\in N$, then $C\nsim D'$. If
$C\notin N$, then $C\thicksim C'= (0, \cdots, 0, c_j, 0 \cdots,
0)$ for some $j\neq i$ by (\ref{eq2}), so that  $C'$ and $D'$ are
different vertices in $N\cup A$, it follows that $C\nsim D$.\\
Case~4. $C, D\in T_1$. In this case, let $i$ be the least integer
for which $c_i\notin S_i$ and let $j$ be the least integer for
which $d_j\notin S_j$; then $C\thicksim C'=(0, \cdots, 0,
\tilde{c_i}, 0 \cdots, 0)$ and $D\thicksim D'=(0, \cdots, 0,
\tilde{d_j}, \cdots, 0)$ by (\ref{eq3}). If $i=j$, then
$c_id_i=0$, contradicts to the assumption (ii). Therefore, $i\neq
j$. It follows that $C'$ and $D'$ are different vertices in $N\cup
A$, $C'\nsim D'$. Thus, $C\nsim D$. \\ Case~5. $C\in T_1$ and
$D\in T_2$. In this case, let $i$ be the least integer for which
$c_i\notin S_i$ and let $j$ be the least integer for which
$d_j\notin Z(R_j)$; then $C\thicksim C'=(0, \cdots, 0,
\tilde{c_i}, 0 \cdots, 0)$  by (\ref{eq3}) and $D\thicksim D'=(0,
\cdots, 0, b_j, \cdots, 0)$ by (\ref{eq4}). If $i=j$, then
$c_id_i=0$, contradicts to the assumption (ii). Therefore, $i\neq
j$. It follows that $C'$ and $D'$ are different vertices in $N\cup
A$, $C'\nsim D'$. Thus, $C\nsim D$. \\ Case~6. $C, D\in T_2$. In
this case, let $i$ be the least integer for which $c_i\notin
Z(R_i)$ and let $j$ be the least integer for which $d_j\notin
Z(R_j)$; then $C\thicksim C'=(0, \cdots, 0, b_i, 0 \cdots, 0)$ and
$D\thicksim D'=(0, \cdots, 0, b_j, \cdots, 0)$ by (\ref{eq4}). If
$i=j$, then $c_id_i=0$, contradicts to the assumption (ii).
Therefore, $i\neq j$. It follows that $C'$ and $D'$ are different
vertices in $N\cup A$, $C'\nsim D'$. Thus, $C\nsim D$.
\end{proof}

\begin{coll}
\label{maincoll1} Let $R\cong \mathbb{Z}^n\times
\mathbb{Z}_{p_1^{r_1}}\times \cdots \times
\mathbb{Z}_{p_k^{r_k}}$; then
$\chi(\Gamma_0(R))=\omega(\Gamma_0(R))=\dprod_{i=1}^k
p_i^{\lfloor{r_i/2}\rfloor}+t+n$, where $\lfloor{n}\rfloor$ is the
largest integer that is smaller than or equal to $n$ and
$t=|\{i~|~r_i~is~odd\}|$.
\end{coll}
\begin{proof}
By Lemma~\ref{lem2}, Lemma~\ref{lem3} and Corollary~\ref{coll1},
we may assume that $R\cong \mathbb{Z}_{p_1^{r_1}}\times \cdots
\times \mathbb{Z}_{p_k^{r_k}}$. Let $S_i=\{p_i^{r_i/2}n~|~n\in
\mathbb{Z}\}$ if $r_i$ is even and $S_i=\{p^{(r_i+1)/2}n~|~n\in
\mathbb{Z}\}\cup \{p^{(r_i-1)/2}\}$ if $r_i$ is odd; then $S_i$
satisfies (i) and (ii) of Theorem~\ref{maintheo1}. Let $N_i=\{a\in
S_i~|~a^2=0\}$; then by Lemma~\ref{lem1} and
Corollary~\ref{coll1}, $|S_i-N_i|=1$ and
$|N_i|=p^{\lfloor{r_i/2}\rfloor}$ if $r$ is odd. Further,
$|S_i-N_i|=0$ and $|N_i|=p^{\lfloor{r_i/2}\rfloor}$ if $r$ is
even. Thus, the assertion follows.
\end{proof}

Since $\Gamma(R)$ is a subgraph of $\Gamma_0(R)$ and
$\Gamma_0(R)-\{0\}$ is a disjoint union of $\Gamma(R)$ and a
totally disconnected graph. We have the following.
\begin{coll}
\label{maincoll2} Let $R\cong \mathbb{Z}_{p_1^{r_1}}\times \cdots
\times \mathbb{Z}_{p_k^{r_k}}$; then
$\chi(\Gamma(R))=\omega(\Gamma(R))$.
\end{coll}

\begin{remark}
\emph{(i) In \cite[Theorem~3.2]{an}, Anderson and Naseer also
provide some sufficient conditions for the product of coloring
rings to be chromatic. However, \cite[Theorem~3.2]{an} can not
derive Theorem~\ref{maintheo1}.  For example, if $R\cong
\mathbb{Z}^n$, then Theorem~\ref{maintheo1} can be used to obtain
that
$R$ is chromatic. However, \cite[Theorem~3.2]{an} can not. \\
(ii) Theorem~\ref{maintheo1} can be used to derive
\cite[Corollary~3.3]{an}. \\ (iii) In \cite{b}, it is shown that
if $R\cong \mathbb{Z}_n$, then $R$ is chromatic.
Corollary~\ref{maincoll1} generalize this result.}
\end{remark}

\section{Coloring of $\overline{\Gamma(R)}$}

Throughout, let $R\cong \mathbb{Z}_{p_1^{r_1}}\times \cdots \times
\mathbb{Z}_{p_k^{r_k}}$. Let $\Gamma(R)$ be the graph associated
to $R$ defined in the introduction. For convenience, we fix some
notion from \cite{cgh}.

\begin{defi}
\label{defi1} Let $R\cong \mathbb{Z}_{p^r}$, where $p$ is a prime
number and $r$ be a positive integer. Let $a, b$ be two distinct
vertices of $\Gamma(R)$; then $a$ and $b$ are associated, denoted
by $a\thicksim b$, if there is a unit $u$ in $R$ such that $a=ub$.
The associate class of $a$, denoted by $A_a$, is the set $\{b\in
\Gamma(R)~|~b\thicksim a\}$. Thus, $\Gamma(R)=\sqcup_{i=1}^r
A_{p^i}$, $1=|A_{p^r}|<|A_{p^{r-1}}|< \cdots < |A_{p^1}|$ and
$|A_{p^i}|=p^{r-i}-p^{r-i-1}$ for $i=0, 1, \dots, r-1$.
\end{defi}

\begin{lemma}
\label{lembar} Let $R\cong \mathbb{Z}_{p^r}$, where $p$ is a prime
number and $r\geq 2$ be a positive integer. Then
$\chi(\overline{\Gamma(R)})=\omega(\overline{\Gamma(R)})=p^{r-1}-p^{r-(r/2)}$
if $r$ is even and
$\chi(\overline{\Gamma(R)})=\omega(\overline{\Gamma(R)})=p^{r-1}-p^{(r-1)/2}$
if $r$ is odd.
\end{lemma}
\begin{proof}
If $r=2$, then $\overline{\Gamma(R)}$ is discrete, so that
$\chi(\overline{\Gamma(R)})=\omega(\overline{\Gamma(R)})=0$,
therefore the assertion holds. Thus, we may assume that $r\geq
3$.\\ Suppose that $r$ is even. Let $C=\sqcup_{i=1}^{r/2-1}
A_{p^i}$; then $|C|=p^{r-1}-p^{r-(r/2)}$ by Definition~\ref{defi1}
and $C$ is a clique of $\overline{\Gamma(R)}$. Let $n=|C|$; then
we require $n$ colors to paint $C$. Since
$|\overline{\Gamma(R)}-C|=p^{r-(r/2)}-1\leq n$, we can use part of
the colors of $C$ to color $\overline{\Gamma(R)}-C$. Hence
$\chi(\overline{\Gamma(R)})\leq |C|\leq
\omega(\overline{\Gamma(R)})$.\\ Suppose that $r$ is odd. Let
$C=\sqcup_{i=1}^{(r-1)/2} A_{p^i}$; then $|C|=p^{r-1}-p^{(r-1)/2}$
by Definition~\ref{defi1} and $C$ is a clique of
$\overline{\Gamma(R)}$. Let $n=|C|$; then we require $n$ colors to
paint $C$. Since $|\overline{\Gamma(R)}-C|=p^{(r-1)/2}-1\leq n$,
we can use part of the colors of $C$ to color
$\overline{\Gamma(R)}-C$. Hence $\chi(\overline{\Gamma(R)})\leq
|C|\leq \omega(\overline{\Gamma(R)})$.
\end{proof}

\begin{example}
\label{exam2} Let $R\cong \mathbb{Z}_8\times \mathbb{Z}_{16}$;
then $\chi(\overline{\Gamma(R)})=\omega(\overline{\Gamma(R)})=76$.
\end{example}

\begin{proof}
Let $$\ba{rl} C= & \{0, 2, 4, 6\}\times \{1, 3, 5, 7, 9, 11, 13,
15\}\\ & \sqcup (\mathbb{Z}_8\times \{2, 6, 10, 14\})\\ &
\sqcup(\{1, 3, 5, 7, 2, 6\}\times \{4, 12\} ) \ea;
$$ then $C$ is a clique in $\overline{\Gamma(R)}$ as $C$ is an
independent set in  $\Gamma(R)$. Moreover, $C$ is maximal. Observe
that $|\Gamma(R)|=8\cdot 16-4\cdot 8-1=95$ and $|C|=76$. So,
$|\overline{\Gamma(R)}-C|=19<|C|$. Now, we require $76$ colors to
paint $C$ and then we can use part of the colors of $C$ to color
$\overline{\Gamma(R)}-C$. Hence, $\chi(\overline{\Gamma(R)})\leq
|C|\leq \omega(\overline{\Gamma(R)})$.
\end{proof}

\begin{theo}
\label{maintheo2} Let $R\cong \mathbb{Z}_{p_1^{r_1}}\times \cdots
\times \mathbb{Z}_{p_k^{r_k}}$, where $p_i$ is a prime number for
every $i$. Then
$\chi(\overline{\Gamma(R)})=\omega(\overline{\Gamma(R)})$.
\end{theo}
\begin{proof}
By Lemma~\ref{lembar}, we may assume that $k\geq 2$. Let
$$T=\{(t_1, \dots, t_k)\in \mathbb{N}^k~|~0\leq t_i\leq
r_i~\forall i\}-\{(0, \dots, 0), (r_1, \dots, r_k)\}.$$ For any
$k$-tuple $(t_1, \dots, t_k)\in T$, let $$A(t_1, \dots,
t_k)=A_{p_1^{t_1}}\times \cdots \times A_{p_k^{t_k}}.$$ Here,
$A_{p_i^{t_i}}$ is the associated class of $p_i^{t_i}$ in
$\mathbb{Z}_{p_i^{r_i}}$. It follows that $$|A(t_1, \dots,
t_k)|=\dprod_{i=1}^k |A_{p_i^{t_i}}|
$$ and
$$\Gamma(R)=\sqcup_{(t_1, \dots, t_k)\in T} A(t_1, \dots, t_k).$$ For any
$k$-tuple $(t_1, \dots, t_k)\in T$, define $$\tilde{A}(t_1, \dots,
t_k)=A_{p_1^{r_1-t_1}}\times \cdots \times
A_{p_k^{r_k-t_k}}=A(r_1-t_1, \dots, r_k-t_k).$$ Observe that if
$r_i$ is even for every $i$, then $\tilde{A}(r_1/2, \dots,
r_k/2)=A(r_1/2, \dots, r_k/2)$. If this is the case, we use $A$
for $A(r_1/2, \dots, r_k/2)$. \par Consider the two following
subsets of $T$: $$T_1=\{(t_1, \dots, t_k)\in T~|~|A(t_1, \dots,
t_k)|> |\tilde{A}(t_1, \dots, t_k)|\}$$ and
$$T_2=\{(t_1, \dots, t_k)\in T~|~|A(t_1, \dots, t_k)|= |\tilde{A}(t_1,
\dots, t_k)|~and~ there ~is~an~i~such~that~t_i\neq r_i/2\}$$ Let
$T_0$ be the subset of $T$ consists of every $k$-tuple of $T_1$
and exactly one $k$-tuple of $\{(t_1, \dots, t_k), (r_1-t_1,
\dots, r_k-t_k)\}$ from $T_2$. Notice that $(r_1/2, \dots,
r_k/2)\notin T_0$ if $r_i$ is even for every $i$. Notice also that
if $(t_1, \dots, t_k)\in T_0$, then there is an integer $i$ such
that $t_i<r_i/2$. For if not, then $t_i\geq r_i/2$ for every $i$.
Therefore $|A_{p_i^{t_i}}|\leq |A_{p_i^{r_i-t_i}}|$ for every $i$.
If $t_j\neq r_j/2$, then $t_j>r_j/2$ for some $j$, so that
$|A_{p_j^{t_j}}|< |A_{p_j^{r_j-t_j}}|$, it follows that $|A(t_1,
\dots, t_k)|<|\tilde{A}(t_1, \dots, t_k)|$, a contradiction. Thus,
$r_i$ is even  and $t_i=r_i/2$ for every $i$, contradiction again.
\par Let $V$ be a vertex in $A$ if $r_i$ is even for every $i$. Further, let
$$C=\sqcup_{(t_1, \dots, t_k)\in T_0} A(t_1, \dots, t_k)\sqcup \{V\}$$ if $r_i$ is even for every $i$
and $$C=\sqcup_{(t_1, \dots, t_k)\in T_0} A(t_1, \dots, t_k)$$ if
$r_i$ is odd for some $i$.\par In the following, we will show that
$C$ is a clique in $\overline{\Gamma(R)}$. For this, we first
examine that $A(t_1, \dots, t_k)$ is a clique in
$\overline{\Gamma(R)}$ if $(t_1, \dots, t_k)\in T_0$. Or
equivalently, to show that $A(t_1, \dots, t_k)$ is an independent
set in $\Gamma(R)$. Suppose not. Then there are vertices $(a_1,
\dots, a_k), (b_1, \dots, b_k)\in A(t_1, \dots, t_k)$ such that
$a_ib_i=0$ for every $i$. However, this implies that $t_i\geq
r_i/2$ for every $i$, contradicts to the fact that there is an
integer $i$ such
that $t_i<r_i/2$ if $(t_1, \dots, t_k)\in T_0$. \\
We next show that if $(c_1, \dots, c_k)=C\in A(t_1, \dots, t_k)$
and $(d_1, \dots, d_k)=D\in A(s_1, \dots, s_k)$, then $(c_1,
\dots, c_k)$ and $(d_1, \dots, d_k)$ are adjacent in
$\overline{\Gamma(R)}$, where $(t_1, \dots, t_k)$ and $(s_1,
\dots, s_k)$ are different $k$-tuples in $T_0$. Or equivalently,
to show that $(c_1, \dots, c_k)$ is not adjacent to $(d_1, \dots,
d_k)$ in $\Gamma(R)$. Suppose not. Then $c_id_i=0$ for every $i$.
This implies that $t_i+s_i\geq r_i$ for every $i$. Since $A(s_1,
\dots, s_k)\neq \tilde{A}(t_1, \dots, t_k)$, the set
$\{i~|~t_i+s_i>r_i\}$ is nonempty. We may assume that
$t_1+s_1>r_1$. Therefore, $t_1\geq 1$ and $s_1\geq 1$. For, if
$t_1=0$ or $s_1£á=0$, then $c_1$ or $d_1$ is a unit, so that
$s_1=r_1$ or $t_1=r_1$, it follows that $t_1+s_1=r_1$, a
contradiction. Thus,
$$|A(t_1, t_2, \dots, t_k)|<|A(t_1-1, t_2, \dots, t_k)|\leq
|A(r_1-s_1, \dots, r_k-s_k)|=|\tilde{A}(s_1, \dots, s_k)|$$ and
$$|A(s_1, s_2, \dots, s_k)|<|A(s_1-1, s_2, \dots, s_k)|\leq
|A(r_1-t_1, \dots, r_k-t_k)|=|\tilde{A}(t_1, \dots, t_k)|.$$ Thus,
$|A(t_1, \dots, t_k)|< |A(s_1, \dots, s_k)|$ and $|A(s_1,
\dots, s_k)|<|A(t_1, \dots, t_k)|$, a contradiction. \\
Now, if $r_i$ is even for every $i$, then $A$ is a clique in
$\Gamma(R)$, so that $V$ is not adjacent to every vertex of
$A-\{V\}$ and $A-\{V\}$ is an independent set in
$\overline{\Gamma(R)}$. Moreover, $V$ is adjacent to every vertex
of $A(t_1, \dots, t_k)$ in $\overline{\Gamma(R)}$ if $(t_1, \dots,
t_k)\in T_0$ as there is an integer $i$ such that $t_i<r_i/2$ and
then $V$ is not adjacent to every vertex of $A(t_1, \dots, t_k)$
in $\Gamma(R)$. \par  From the above, we see that $C$ is a clique
in $\overline{\Gamma(R)}$. Let $n=|C|$; then we require $n$
different colors to paint $C$. We then use the colors of $A(t_1,
\dots, t_k)$ ($(t_1, \dots, t_k)\in T_0$) to color $\tilde{A}(t_1,
\dots, t_k)$ and the color of $V$ to color every vertex of
$A-\{V\}$ if $r_i$ is even for every $i$. In this way, it is easy
to see that every two adjacent vertices in $\overline{\Gamma(R)}$
have different colors. Hence, we conclude that
$\chi(\overline{\Gamma(R)})\leq n\leq
\omega(\overline{\Gamma(R)})$.
\end{proof}

An easy consequence of Theorem~\ref{maintheo2} is that if $R\cong
\mathbb{Z}_{p_1}\times \cdots \times \mathbb{Z}_{p_k}$, then
$\chi(\overline{\Gamma(R)})=\omega(\overline{\Gamma(R)})$.
However, from the proof of Theorem~\ref{maintheo2}, one can
replace each $\mathbb{Z}_{p_i}$ by any finite field.
\begin{coll}
Let $R\cong F_1\times \cdots \times F_k$, where $F_i$ is a finite
fields for every $i$. Then
$\chi(\overline{\Gamma(R)})=\omega(\overline{\Gamma(R)})$.
\end{coll}

~\\
Department of Mathematics, \\
National Chung Cheng University,\\
Chiayi 621, Taiwan \\
hjwang@math.ccu.edu.tw

\end{document}